\documentclass[12pt,a4paper]{article}
\usepackage{amssymb,amsmath}
\usepackage[T1]{fontenc}
\usepackage{theorem}
{%\theoremstyle{marginbreak}
{\theorembodyfont{} %kroj pisma w tw
\newtheorem{thm}{Theorem}[section]
\newtheorem{cor}[thm]{Corollary}

\newtheorem{prop}[thm]{Proposition}
{\theorembodyfont{\normalfont}
\newtheorem{defin}[thm]{Definition}
\newtheorem{rem}[thm]{Remark}
\newtheorem{exa}[thm]{Example}
}} %dziala lokalnie

\title{Conformal fields and the stability of leaves with constant higher order mean
curvature}
\author{Krzysztof Andrzejewski and Pawe\l\, G. Walczak
}
\date{}
\begin{document}
\def\mR{\mathbb R}
\def\calf{\mathcal F}
\def\cala{\mathcal A}
\def\ric{\operatorname{Ric}}
\def\tr{\operatorname{Tr}}
\newcommand\brac[1]{\langle{#1}|}
\newcommand\ket[1]{|{#1}\rangle}
\newcommand\bracket[2]{\langle{#1},{#2}\rangle}
\newcommand\codim{\operatorname{codim}}
\newcommand\divergence{\operatorname{div}}
\newcommand\vol{\operatorname{vol}}

\maketitle

\begin{abstract}
In this paper, we study submanifolds  with constant $r$th mean curvature $S_r$.
We investigate,  the stability of such submanifolds in the case when they  are leaves of a codimension
 one foliation. We also generalize recent results by Barros - Sousa and Al\'{i}as - Colares,  concerning conformal fields, to an arbitrary manifold.
Using this   we show that normal component of a Killing field is a $r$th Jacobi field of a submanifold
with $S_{r+1}$ constant. Finally, we study  relations between $r$th Jacobi fields and vector fields preserving
a foliation.
\end{abstract}

\section{Introduction}
Let $L$ be a  submanifold of $M$ with a unit normal vector field $N$ and  constant $(r+1)$th mean curvature $S_{r+1}$.
If  $M$ is the manifold with  constant sectional curvature (Einstein manifold  for $r=1$) then $L$ is
characterized by a variational problem (see, among the others, \cite{ace,bc,bs2,e}). Therefore, there is a natural
question about the stability of $L$. In this paper, we give some criteria for  the stability
of such submanifolds in the case when they  are leaves of a codimension
one foliation (Theorem \ref{t:1}). Next, for a conformal vector field $U$ we obtain  the formula for  $L_r(f)$,
where $f=\langle U,N\rangle$, in the case of an arbitrary manifold (Theorem \ref{t:2}).
Using this,  we show that normal component of a Killing field is a $r$th Jacobi field of a submanifold
with $S_{r+1}$ constant  (Proposition \ref{p:2}). Finally,
we investigate  relations between $r$th Jacobi fields and vector fields preserving
a foliation - Section \ref{s:3}.
\par
Throughout the paper everything
(manifolds, foliations, metrics, etc.) is assumed to be
$C^{\infty}$-differentiable and oriented.
For simplicity, we usually work with $S_r$ instead of its
normalized counterpart $H_r$ (see  Remark \ref{r1}).
Repeated indices denote summation over their range.
\section{Preliminaries}
\label{s:0}
Let $M$ be a $(n+1)$-dimensional
Riemannian  manifold, $L$ be a codimension one  submanifold
of $M$ and $\langle \cdot ,\cdot\rangle$  represent a  metric on $M$.
Assume that both $M$ and $L$ are oriented and let $N$ be an  orthogonal  unit vector
field. Let $\overline\nabla$ denote the Levi-Civita connection
of the metric. Then $\overline\nabla $ induces  the  connection $\nabla$  on
the set $\Gamma(L)$ of all vector  fields  tangent  to
$L$. Define the second fundamental form (or, the shape operator) $A$ of
$L $ with respect to $N$ by
\[
A: \Gamma(L) \rightarrow
\Gamma(L), \quad A(X)=-(\overline\nabla_XN )^\top \quad \textrm{for } X\in
\Gamma(L),
\]
where $^\top $ denotes the orthogonal  projection on the
vector bundle tangent to $L$.  Note that  $A$ is a self-adjoint linear
operator and at each  point $p\in L$   has  real eigenvalues $\kappa
_1(p),\ldots,\kappa_n(p)$ (the principal curvatures). Associated to the
shape operator there are $n$ algebraic invariants given by
\[
 S_r(p)=\sigma_r(\kappa_1(p),\ldots,\kappa_n(p)) ,
\]
where $\sigma_r$ for
$r=1,2,\ldots,n$  are the elementary symmetric functions given by
\[
 \sigma _r(x_1,\ldots,x_n) =\sum_{i_1<\cdots<i_r}x_{i_1}\cdots x_{i_r},
 \]
$\sigma_0=1$ and $\sigma_r=0$ for other $r$. Moreover, observe that the
characteristic polynomial of $A$ can be written  in terms of the $S_r$'s as
\[
\det(tI-A)=\sum _{r=0}^{n} (-1)^rS_rt^{n-r}.
\]
The normalized $r$th mean
curvature $H_r$ of $L$ is defined by
\[
 H_r=S_r\dbinom{n}{r}^{-1}.
\]
\begin{rem}
 \label{r1} {\rm Sometimes $H_r$,
instead of $S_r$,  is called $r$th mean curvature.}
 \end{rem}
\par
Now, we  introduce the Newton
transformations $T_r:\Gamma(L)\rightarrow \Gamma(L)$ arising from
the shape operator. They  are defined inductively by
\[
T_0=I,\quad T_r=S_rI-AT_{r-1}, \quad 1\leq r\leq n,
\]
or, equivalently, by
\[
 T_r = S_rI-S_{r-1}A+\cdots+(-1)^{r-1}S_1A^{r-1}+(-1)^rA^r.
\]
Note that, by the
Cayley-Hamilton theorem we have $T_n$=0.  Furthermore, $T_r$ is also
self-adjoint and $A$ together with all the $T_r$'s can be simultaneously
diagonalized; if
$e_1,\ldots,e_n$ are the eigenvectors of $A$ corresponding to the
eigenvalues $\kappa _1(p),\ldots ,\kappa _n(p)$, respectively, then they are
also eigenvectors of $T_r$ corresponding to the eigenvalues $\mu_{i,r}(p)$
of $T_r$, that is $T_r(e_i)=\mu_{i,r}(p)e_i$, where
\[
\mu_{i,r}(p)=\frac{\partial \sigma_{r+1}}{\partial
x_i}(\kappa_1(p),\ldots,\kappa_n(p)).
\]
We say that $T_r$ is definite (semi definite) if $T_r>0$ or $T_r<0$  on $L$
 ($T_r\geq 0$ or $T_r\leq0$ on $L$).
\par The following algebraic properties
of $T_r$ are well known (see, for instance, \cite{ro}) and will be applied
throughout this paper:
\begin{align*}
&\tr(T_r)=(n-r)S_r =c_rH_r,\\
&\tr(AT_r)=(r+1)S_{r+1}=c_rH_{r+1},\\
&\tr(A^2T_r)=S_1S_{r+1}-(r+2)S_{r+2},
\end{align*}
where $c_r=(n-r)\dbinom{n}{r}=(r+1)\dbinom{n}{r+1}$.
\par
Let   $f\in C^\infty(L)$.
Define  operators $L_r,J_r$ as follows:
\[
L_rf=\tr(T_r\circ {\rm Hess} f),
\]
and
\[
J_rf=L_rf+\tr(A^2T_r)f+\tr(\overline R(N)T_r)f,
\]
where $\overline R(N): \Gamma(L)\rightarrow \Gamma(L)$ is given by
\[
\overline R(N)(X)=\overline R(X,N)N, \quad X\in\Gamma(L),
\]
and $\overline R$ being the curvature tensor of $\overline \nabla$.
Then
\[
L_rf=\divergence(T_r{\rm \nabla}f)-\langle\divergence T_r,{\rm \nabla} f\rangle,
\]
where $\divergence  T_r =(\nabla_{e_i}T_r)e_i$,
and we have the following cases (see, among the others, \cite{alm,bc,bs,cr,e}).
\par
For $r=0$ we have $\divergence_L T_0=0$ thus  $L_r=L_0=\Delta$
\[
J_0f=\Delta f+\tr A^2f +\overline\ric(N)f.
\]
If $r=1$ and  $M$  is an Einstein manifold,  then $\divergence T_1=0$ and
\[
J_1f=\divergence(T_1{\rm \nabla}f)+(S_1S_{2}-3S_{3})f+\tr(\overline R(N)T_1)f.
\]
If $M$  is a manifold with constant sectional  curvature $c$, then  for arbitrary $r$
$\divergence (T_r)=0$ and
\[
J_rf=\divergence(T_r{\rm \nabla}f)+(S_1S_{r+1}-(r+2)S_{r+2})f+(n-r)cS_rf.
\]
For these three cases we have the following proposition (e.g. \cite{bc}).
\begin{prop}
\label{p0}
If $L$ is compact without boundary or if $L$ is noncompact and $f\in C_c^\infty(L)$  then
\[
\int _L L_r(f)=0 \quad {\rm and} \quad \int_L fL_r(f)=-\int_L\langle T_r\nabla f,\nabla f\rangle.
\]
\hfill$\square$
\end{prop}
Next, we define
\[
I_r(f,g)=-\int_L fJ_rg,
\]
for $f,g\in C_c^\infty(L)=\{f\in C^\infty(L): f-{\rm\, a\, compactly\, supported\, function}\}$.
\par
Let us  recall that a submanifold is $r$-minimal   ($0\leq r\leq n-1$) if $S_{r+1}=0$.
Let $\calf$ be  a codimension one foliation. We say that  $\calf$ is $r$-minimal  if
any leaf of $\calf$  is a $r$-minimal submanifold of $M$.   A foliation such that
every leaf  has constant $(r+1)$th mean  curvature is called $r$-tense.
\par  Similarly, as for submanifolds, we may define $S_r,H_r,T_r$ for
a foliation $\calf$ (e.g. \cite{aw,aw2}).  In this case, the  functions
$S_r$,   are  smooth  on  the whole $M$ and, for any point
$p\in M$,  $S_r(p)$  coincides with the r-th mean curvature at $p$  of the leaf
$L$ of $\calf$ which passes through  $p$; therefore we will use the same notation for
$r$th mean curvature of foliations and submanifolds.
Finally, recall that a hypersurface $L$  with  $S_{r+1}=constant$,   of a manifold with  constant sectional curvature
(Einstein manifold -- for $r=1$) is a critical point of  the variational problem of minimizing the integral
\[
\cala_r=\int_LF_r(S_1,\ldots,S_r),
\]
for compactly supported volume-preserving variations, see \cite{bc,bs2,e}.
The functions $F_r$ are defined inductively by
\begin{align*}
&F_0=1,\\
&F_1=S_1,\\
&F_r=S_r+\frac{c(n-r+1)}{r-1}F_{r-2}, \quad 2\leq r\leq n-1.
\end{align*}
The second variation formula reads $\cala_r''(0)=(r+1)I_r(f,f)$.
Thus,  we  may introduce the following definition (see discussion in \cite{ace}).
 \begin{defin}We say that a submanifold $L$ with $S_{r+1}=constant$ is $r$-stable if $I_r(f,f)\geq 0$ for all $f\in C_c^*(L)$ or if
$I_r(f,f)\leq 0$ for all $f\in C_c^*(L)$. We say $L$ is $r$-unstable if there exist  functions
$f,g\in C_c^*(L)$ such that $I_r(f,f)<0$  and  $I_r(g,g)>0$;
where
\begin{equation}
\label{e0}
C_c^*=\{f\in C^\infty_c: \int_Lf=0\}.
\end{equation}
\end{defin}
$0$-minimal ($0$-stable) submanifold are simply called minimal (stable).

\section{Stability results}
\label{s:1}
Oshikiri \cite{o1} has showed that  each leaf of  a minimal foliation is stable.
Now, we give a generalization of  this theorem for arbitrary $r>0$.
In order to do this, we will need the following proposition (\cite{aw}, see also \cite{csc}).
\begin{prop}
\label{p1}
Let $M$ be a  Riemannian manifold with a unit
vector field $N$ orthogonal to the foliation $\calf$ of $M$. Then on a  leaf $L$ we have
\begin{align*}
\divergence(T_r\overline\nabla_NN)&=\langle
\divergence T_r,\overline\nabla_NN\rangle-N(S_{r+1})+\nonumber\\
&+\tr(A^2T_r)+\tr(\overline R(N)T_r)+\langle\overline\nabla_NN,T_r\overline\nabla_NN\rangle.
\end{align*}
\hfill$\square$
\end{prop}
%As a vector field $\nabla_NN$ is tangent to $\calf$ everywhere, the restriction of it  to nay leaf
%$L$ is well-defined vector field on $L$ and we also denote it  by $\nabla_NN$.
If $M$ is a manifold without boundary, then we have the following theorem.
\begin{thm}
\label{t:1}
Let $M$ be an Einstein (a constant sectional curvature)  manifold   for $r=1$ ($r>1$) and
 $\calf$ be a  codimension one  $r$-tense foliation of $M$.
If on a leaf $L$, either $T_r\geq 0$ and $N(S_{r+1})\leq 0$  or $T_r\leq 0$ and $N(S_{r+1})\geq 0$,
then  $L$ is $r$-stable.
\end{thm}
{\it Proof.}
In our case from Proposition \ref{p1}  we get
\begin{align*}
\divergence(T_r\overline\nabla_NN)&=-N(S_{r+1})+\tr(A^2T_r)+\tr(\overline R(N)T_r)+
\langle\overline\nabla_NN,T_r\overline\nabla_NN\rangle.
\end{align*}
Thus, for any $f\in C_c^*(L)$ we have
\begin{align}
\label{e1}
&\divergence(f^2T_r\overline\nabla_NN)-(T_r\overline\nabla_NN)(f^2)=
f^2\divergence(T_r\overline\nabla_NN)
\nonumber\\
&=f^2\tr(T_rA^2)+f^2\tr(\overline R(N)T_r)+f^2\langle\overline\nabla_NN,T_r\overline\nabla_NN\rangle
-f^2N(S_{r+1})
\end{align}
Using Proposition \ref{p0},  Eq. (\ref{e1})  and the fact that $T_r$ is selfadjoint,  we have
\begin{align*}
&\int_L\langle T_r(\nabla f+f\overline\nabla_NN),\nabla f+f\overline\nabla_NN\rangle
-f^2N(S_{r+1})\\
&=\int_L\langle T_r\nabla f,T_r\nabla f\rangle+2f\langle T_r\nabla_NN,\nabla f\rangle
+ f^2\langle T_r\overline\nabla_NN,\overline\nabla_NN\rangle-f^2N(S_{r+1})\\
&=\int_L-fL_r(f)+(T_r\overline\nabla_NN)(f^2)+ f^2\langle T_r\overline\nabla_NN,
\overline\nabla_NN\rangle-f^2N(S_{r+1})\\
&=\int_L-fL_r(f)-f^2\tr(T_rA^2)-f^2\tr(R(N)T_r)+\divergence(f^2\overline\nabla_NN)\\
&=I_r(f,f).
\end{align*}
This ends the proof.\hfill$\square$
\medskip
\par
Note that,  during the proof of Theorem \ref{t:1}, we did not use the condition from Eq. (\ref{e0}).
\begin{cor}
Let $M$ be as in Theorem \ref{t:1}. If each leaf of the foliation $\calf$ has the same constant $(r+1)$th
mean curvature (especially equal zero) and $T_r$ is semi definite on $M$ then any leaf of $\calf$ is $r$-stable.
\end{cor}
There are various conditions enforcing (semi) definiteness of the operator $T_r$, see
\cite{asz,cr}. One of them implies the following corollary.
\begin{cor}
Let $M$ be as in Theorem \ref{t:1}  and $\calf$ be a $r$-minimal foliation of $M$. If on a leaf $L$,
$S_r\neq 0$,  then  $L$ $r$-stable.
\end{cor}
\begin{exa}
Let $M=\mR\times L$ be a foliated manifold each leaf of which is given by $\{t\} \times L$
where $L$ has constant negative sectional curvature $c$. We define a metric on $M$  by $\langle,\rangle=
dt^2+\cosh(\sqrt{-c}t)\langle,\rangle_L$. Then $(M,\langle,\rangle)$ has constant sectional curvature $c$
and foliation is $r$-tense. Moreover,  on any leaf $L$, $T_r\geq 0$ and $N(S_{r+1})\leq 0$  or $T_r\leq 0$
and $N(S_{r+1})\geq 0$,  thus  any leaf is $r$-stable.
\end{exa}
\begin{exa}
Let $M=\mR\times L$ be a foliated manifold each leaf of which is given by $\{t\} \times L$,
where $L$  is flat manifold (e.g. $\mR^n$,$T^n$). We define a metric on $M$  by $\langle,\rangle=
dt^2+e^{-2at}\langle,\rangle_L$. Then $(M,\langle,\rangle)$ has constant sectional curvature $-a^2$,
each leaf has the same constant $S_{r+1}$ and $P_r$ is definite;
thus each leaf is $r$-stable.
\end{exa}
Recall that  by a singular foliation of $M$, we mean a foliation  $\calf$
of $M\backslash S$, where $S\subset M$ is a set of Lebesgue measure zero \cite{csc}.
\begin{exa}
Let $\calf$ be a  singular foliation of $\mR^{n+1}$ by the concentric cylinders $S^r(R)\times \mR^{n-r}$,
where $S^r(R)$  denotes the sphere with center $0\in \mR^{n+1}$ and  radius $R > 0$; the singular
set of the foliation is the $(n- r)$-hyperplane $\{0\}\times \mR^{n-r}$ in $\mR^{n+1}$. Then $\calf$
is $r$-minimal foliation and $S_r\neq 0$; consequently,  any leaf is $r$-stable.
\end{exa}
\section{Conformal fields}
\label{s:2}
Let $U$ be a conformal vector field on a manifold $M$ and $f=\langle U,N\rangle$.
Recently,  Barros - Sousa and Al\'{i}as - Colares  \cite{ac,bc} have obtained an expression of $L_rf$
when  $M$ is either a manifold with  constant sectional curvature or generalized a Robertson-Walker spacetime.
Now, we generalize these results  to the case of arbitrary manifolds and obtain some other consequences.
\begin{thm}
\label{t:2}
Let $L$ be a submanifold (not necessary a leaf) of an arbitrary manifold $M$ with the unit normal
vector field $N$. If  $U$ is a conformal vector field  on $M$ and  $f=\langle U,N\rangle$,
 then
\[
J_rf=-U^\top(S_{r+1})-(r+1)kS_{r+1}-N(k)(n-r)S_r,
\]
equivalently
\[
L_rf=-\langle U,\nabla S_{r+1}\rangle -f\tr(A^2T_r)-f\tr(\overline R(N)T_r)-k\tr(AT_r)-N(k)\tr(T_r),
\]
where $2k$ is  the conformal factor of $U$.
\end{thm}
{\it Proof.}
Let  $X\in\Gamma(T(L))$. Since $U$ is a conformal field,   we have
\begin{align*}
&\langle \nabla f,X\rangle=X(f)=\langle\overline \nabla_XU,N\rangle+\langle U,\overline\nabla_XN\rangle\\
&=-\langle X,\overline\nabla_NU\rangle+\langle U^\top,\overline\nabla_XN\rangle\\
&=-\langle X,\overline\nabla_NU \rangle-\langle U^\top,AX\rangle=
-\langle X,(\overline\nabla_NU)^\top+AU^\top\rangle.
\end{align*}
Thus, we get
\begin{equation}
 \nabla f =-((\overline\nabla_NU)^\top+AU^\top).
\end{equation}
Let $p\in L$ be  an arbitrary point and $\{e_i\}_{i=1}^n$ a local  orthonormal frame such that
 $T_r(e_i(p))=\mu_{i,r}e_i(p)$.
By definition of $L_r$ we have
\[
(L_rf)(p)=\langle\nabla_{e_i}(\nabla f),T_re_i\rangle(p),
\]
where as everywhere, repeated indices denote summation.
\par
Thus   at the point $p$  we obtain
\begin{align*}
&\langle \nabla_{e_i}(\overline\nabla_NU)^\top,T_re_i\rangle =
\langle \overline\nabla_{e_i}(\overline\nabla_NU)^\top,T_re_i\rangle \\
&=\langle \overline\nabla_{e_i}\overline\nabla_NU,T_re_i\rangle-\langle\overline\nabla_NU,N\rangle
\langle \overline\nabla_{e_i}N,T_re_i\rangle \\
&=\langle\overline R(e_i,N)U,T_re_i\rangle+k\tr(AT_r)
+\langle\overline\nabla_N\overline\nabla_{e_i}U,T_re_i\rangle+
\langle\overline\nabla_{[e_i,N]}U,T_re_i\rangle\\
&=\langle\overline R(e_i,N)U,T_re_i\rangle+k\tr(AT_r)+
\langle\overline\nabla_N\overline\nabla_{e_i}U,T_re_i\rangle\\
&+\langle\overline\nabla_{\overline\nabla_{e_i}N}U,T_re_i\rangle-
\langle\overline\nabla_{\overline\nabla_Ne_i}U,T_re_i\rangle\\
&=\langle\overline R(e_i,N)U,T_re_i\rangle+k\tr(AT_r)+
\langle\overline\nabla_{\overline\nabla_{e_i}N}U,T_re_i\rangle\\
&+\mu_{i,r}(\langle\overline\nabla_N\overline\nabla_{e_i}U,e_i\rangle-
\langle\overline\nabla_{\overline\nabla_Ne_i}U,e_i\rangle).
\end{align*}
Since, for a fixed $i$ we have  $\langle\overline\nabla_{e_i}U,e_i\rangle=k$ and
$\langle\overline\nabla_Ne_i,e_i\rangle$=0 thus
\[
\langle\overline\nabla_N\overline\nabla_{e_i}U,e_i\rangle=
-\langle\overline\nabla_{e_i}U,\overline\nabla_Ne_i\rangle+N(k)
=\langle\overline\nabla_{\overline\nabla_Ne_i}U,e_i\rangle+N(k).
\]
Consequently, at $p$ we have
\begin{align}
\label{e2}
&\langle \nabla_{e_i}(\overline\nabla_NU)^\top,T_re_i\rangle\nonumber\\
&=\langle\overline R(e_i,N)U,T_re_i\rangle+
\langle\overline\nabla_{\overline\nabla_{e_i}N}U,T_re_i\rangle
+k\tr(AT_r)+N(k)\tr(T_r)\nonumber\\
&=\langle\overline R(T_re_i,N)U,e_i\rangle-
\langle Ae_i,e_j\rangle \langle(\overline\nabla_{e_j}U)^\top,T_re_i\rangle
+k\tr(AT_r)+N(k)\tr(T_r)\nonumber\\
&=\langle\overline R(e_i,U)N,T_re_i\rangle-\tr(AT_r(\overline\nabla U)^\top)
+k\tr(AT_r)+N(k)\tr(T_r)
\end{align}
On the other hand, from the Codazi equation, we obtain
\begin{align*}
&\langle\overline R(e_i,U^\top)N,T_re_i\rangle=\langle(\nabla_{U^\top}A)e_i,T_re_i\rangle-
\langle(\nabla_{e_i}A)U^\top,T_re_i\rangle\\
&=\langle (T_r\nabla_{U^\top}A)e_i,e_i\rangle-\langle\nabla_{e_i}(AU^\top),T_re_i\rangle+
\langle A(\nabla_{e_i}U^\top),T_re_i\rangle\\
&= \tr(T_r\nabla_{U^\top}A)-\langle\nabla_{e_i}(AU^\top),T_re_i\rangle+
\langle \overline\nabla_{e_i}U^\top,AT_re_i\rangle\\
&=U^\top(S_{r+1})-\langle\nabla_{e_i}(AU^\top),T_re_i\rangle
-\langle \overline\nabla_{e_i}(fN),AT_re_i\rangle+\langle \overline\nabla_{e_i}U,AT_re_i\rangle\\
&=U^\top(S_{r+1})-\langle\nabla_{e_i}(AU^\top),T_re_i\rangle+f\langle Ae_i,AT_re_i\rangle
+\langle \overline\nabla_{e_i}U,AT_re_i\rangle\\
&=U^\top(S_{r+1})-\langle\nabla_{e_i}(AU^\top),T_re_i\rangle+f\tr(A^2T_r)
+\tr(T_rA(\overline\nabla U)^\top).
\end{align*}
Thus
\begin{align}
\label{e3}
\langle\nabla_{e_i}(AU^\top),T_re_i\rangle&=
-\langle\overline R(e_i,U^\top)N,T_re_i\rangle+U^\top(S_{r+1})\nonumber\\
&+f\tr(A^2T_r)
+\tr(T_rA(\overline\nabla U)^\top).
\end{align}
Since   $AT_r=T_rA$, we have
\begin{equation}
\label{e4}
\tr(T_rA(\overline\nabla U)^\top)=\tr(AT_r(\overline\nabla U)^\top).
\end{equation}
Finally,  from Eq. (\ref{e2}),(\ref{e3}) and (\ref{e4})  we get at the point $p$
\begin{align*}
L_rf=&\langle\nabla_{e_i}(\nabla f),T_re_i\rangle
=-\langle\overline R(e_i,U)N,T_re_i\rangle+\langle\overline R(e_i,U^\top)N,T_re_i\rangle\\
-&U^\top(S_{r+1})-f\tr(A^2T_r)-k\tr(AT_r)-N(k)\tr(T_r)\\
=&-f\tr(\overline R(N)T_r)-f\tr(A^2T_r)-U^\top(S_{r+1})
-k\tr(AT_r)-N(k)\tr(T_r).
\end{align*}
Since $p$ is arbitrary,  the assertion follows.\hfill$\square$
\begin{cor}
\label{c:1}
When $U$  is a Killing field we get
\[
J_r(f)=-U^\top(S_{r+1})=-\langle\nabla S_{r+1},U\rangle.
\]
\end{cor}
For  further applications see Proposition \ref{p:2} and Corollary \ref{c:2}.
\section{Jacobi fields}
\label{s:3}
Let $M$ be an arbitrary manifold and $L$ be a submanifold of  $M$ with a unit orthogonal field $N$.
Then the operator $J_r$ induct a new mapping (denotes also $J_r$)  $J_r:\Gamma(T(L)^\bot)\rightarrow
\Gamma(T(L)^\bot)$
 as follows
\[
J_r(fN)=J_r(f)N.
\]
\begin{defin}
We say that $V\in \Gamma(T(L)^\bot)$ is a $r$th Jocobi  field of $L$ if $J_r(V)=0$.
We say that $V\in \Gamma(T(\calf)^\bot)$ is a $r$th Jacobi field of $\calf$ if is a $r$th  Jacobi field
for any leaf $L$ of $\calf$.
\end{defin}
\begin{prop}
\label{p:2}
Let $L$ be a submanifold of   an arbitrary Riemannian manifold $M$, such that $S_{r+1}$   is constant
on $L$, then the normal component $U^\bot$ of a Killing vector field $U$ is a $r$th Jacobi vector field.
\end{prop}
{\it Proof.} The proof follows immediately from Corollary \ref{c:1}.\hfill$\square$
\begin{thm}
\label{t:3}
Let $M$ be an arbitrary Riemannian manifold and $\calf$ be a foliation of $M$ whose leaves have the same
constant $(r+1)$th mean curvature (e.g. zero). If $V\in\Gamma(TM)$  preserves
$\calf$ (i.e. maps leaves onto leaves)  then $V^\bot=fN$ is a $r$th Jacobi field of $\calf$.
\end{thm}
{\it Proof.} Since $V$ is foliation preserving,  $[V,\Gamma(T(\calf))]\subset
\Gamma(T(\calf))$ so $\nabla f+f \overline\nabla_NN=0$ on any leaf $L$.
Using this and Proposition \ref{p1}, we get
\begin{align*}
J_rf=&L_rf+f\tr(A^2T_r)+f\tr(\overline R(N)T_r)\\
=&\divergence(T_r(\nabla f))-\langle\divergence T_r,\nabla f\rangle+f\tr(A^2T_r)+f\tr(\overline R(N)T_r)\\
=&\divergence(T_r(\nabla f))+f\divergence(T_r\overline\nabla_NN)-\langle\divergence T_r,
\nabla f+f\overline\nabla_NN\rangle\\
-&f\langle\overline\nabla_NN,T_r\overline\nabla_NN\rangle\\
=&\divergence(T_r(\nabla f+f\overline\nabla_NN))-\langle\divergence T_r,
\nabla f+f\overline\nabla_NN\rangle\\
-&\langle \nabla f+f\overline\nabla_NN,T_r\overline\nabla_NN\rangle
=0.
\end{align*}
\hfill$\square$
\begin{exa}
\label{ex:1}
Let $M=\mR\times \mR^n$ ($M=\mR\times T^n$)  be a foliated manifold  which leaves are $\{t\}\times \mR^n$.
For  functions $\phi_1,\ldots,\phi_n:\mR\rightarrow \mR$ we may define a metric $\langle,\rangle$ on $M$
\[
\langle,\rangle=dt^2+e^{-2\int\phi_i(t)dt}(dx^i)^2.
\]
Then, $S_{r+1}=\sigma_{r+1}(\phi_1,\ldots,\phi_n)$.  So we have a lot of metrics
such that $S_{r+1}$ is constant on $M$. Then, a vector filed $V=f(t)\frac{\partial}{ \partial t} $
is  $\calf$ foliation preserving and consequently $V$ is $r$th Jacobi field. Note that,
$V$ is not a Killing field (in general)  and we could not  use Proposition \ref{p:2}.
\end{exa}
\begin{prop}
\label{p3}
Let  $\calf$ be  a foliation of a Riemannian manifold $M$ whose leaves are closed and have the same constant $(r+1)$th mean curvature.
If on any leaf $L$ operator  $T_r$ is definite  and $\divergence(T_r)=0$, then any $r$th Jacobi fields of $\calf$ preserves $\calf$.
\end{prop}
{\it Proof.}
If  $V=fN$ is a $r$th Jacobi field, then $J_r(V)=0$  thus $I_r(f,f)=0$ on each leaf $L$.
On the other, as in the  proof of Theorem \ref{t:1}, we get
\begin{align*}
&\int_L\langle T_r(\nabla f+f\overline\nabla_NN),\nabla f+f\overline\nabla_NN\rangle=I_r(f,f).
\end{align*}
Thus $V$ is foliation preserving. \hfill$\square$
\medskip
\par
For example, if $M=\mR\times T^n$ and  the metric and foliation  are as in Example \ref{ex:1},
then $\nabla_F A=0$ for any $F\in\Gamma(T(L))$. Consequently $\divergence(T_r)=0$ on
any leaf (although $M$ need not be a manifold with constant sectional  curvature).
\begin{cor}
\label{c:2}
Under  the assumptions of Proposition \ref{p3}, if $U$ is a Killing field on $M$
then $U$ preserves $\calf$.
\end{cor}

{\it Krzysztof Andrzejewski} (corresponding author) \\
Institute of Mathematics, Polish Academy of Sciences\\
ul. \'Sniadeckich 8, 00-956 Warszawa, Poland\\
{\it and}\\
Department of Theoretical Physics II,
University  of \L \' od\'z\\
ul. Pomorska 149/153, 90 - 236  \L \' od\'z, Poland.\\
e-mail: k-andrzejewski@uni.lodz.pl
\vspace{0.5cm}
\\
{\it Pawe\l\, G. Walczak }\\
Institute of Mathematics, Polish Academy of Sciences\\
ul. \'Sniadeckich 8, 00-956 Warszawa, Poland\\
{\it and}\\
Faculty of Mathematics and Informatics, University of \L \' od\'z\\
ul. Banacha 22, 90-238 \L \' od\'z, Poland\\
 e-mail: pawelwal@math.uni.lodz.pl
\end{document}